\RequirePackage{fix-cm}
\documentclass[smallextended]{svjour3}       
\smartqed  
\usepackage{graphicx}
\usepackage{amsmath}
\usepackage{amsfonts}
\usepackage{amssymb}
\usepackage{array,longtable}
\usepackage{amscd}
\usepackage[cp1251]{inputenc}
\usepackage[english]{babel}
\usepackage[usenames]{color}
\usepackage{graphicx}
\graphicspath{{pictures/}}
\DeclareGraphicsExtensions{.pdf,.png,.jpg}
\usepackage[matrix,arrow,curve]{xy}
\begin{document}
\setcounter{tocdepth}{5}
\title{Some bounds for determinants of relatively $D$-stable matrices 
}

\titlerunning{Some bounds for determinants of relatively $D$-stable matrices}        

\author{Olga Y. Kushel
}


\institute{O. Kushel \at
              Shanghai University, \\ Department of Mathematics, \\ Shangda Road 99, \\ 200444 Shanghai, China \\
              Tel.: +86(021)66133292\\
              \email{kushel@mail.ru}}

\date{Received: date / Accepted: date}

\maketitle

\begin{abstract}
In this paper, we study the class of relatively $D$-stable matrices and provide the conditions, sufficient for relative $D$-stability. We generalize the well-known Hadamard inequality, to provide upper bounds for the determinants of relatively $D$-stable and relatively additive $D$-stable matrices. For some classes of $D$-stable matrices, we estimate the sector gap between matrix spectra and the imaginary axis. We apply the developed technique to obtain upper bounds for determinants of some classes of $D$-stable matrices, e.g. diagonally stable, diagonally dominant and matrices with $Q^2$-scalings. 

\end{abstract}

\section{Introduction} Here, as usual, ${\mathcal M}^{n \times n}$ denotes the set of $n \times n$ matrices with real entries, $\sigma({\mathbf A})$ denotes the spectrum of the matrix $\mathbf A$ (i.e. the set of all eigenvalues of $\mathbf A$ defined as zeroes of its characteristic polynomial $f_{\mathbf A}(\lambda):= \det(\lambda{\mathbf I - {\mathbf A}})$).
One of the most prominent and well-known classical results of the theory of positive definite Hermitian matrices is the Hadamard's inequality (see, e.g. \cite{HOJ}, p. 477, Theorem 7.8.1 and \cite{MIR}, p. 417, Theorem 13.5.2).

\begin{theorem} If ${\mathbf A} = \{a_{ij}\}_{i,j = 1}^n$ be a Hermitian positive definite matrix then
\begin{equation} \label{hadeq}
\det{\mathbf A} \leq a_{11}\ldots a_{nn},
\end{equation}
and the equality holds if and only if $\mathbf A$ is diagonal.
\end{theorem}
This celebrated result enjoys a number of proofs, with probably the shortest one is as follows:
\begin{enumerate}
\item[\rm \bf Step 1.] One defines a positive diagonal matrix
 $${\mathbf D} = {\rm diag}\{\frac{1}{a_{11}}, \ \ldots, \ \frac{1}{a_{nn}}\},$$
 and considers ${\mathbf C}:={\mathbf D}{\mathbf A}$ with $c_{ii} = 1$ for all $i = 1, \ \ldots, \ n$. Then \eqref{hadeq} is equivalent to the inequality $\det({\mathbf C}) \leq 1$.
\item[\rm \bf Step 2.] By noticing that ${\mathbf C}={\mathbf D}{\mathbf A}$ is similar to ${\mathbf D}^{\frac{1}{2}}{\mathbf A}{\mathbf D}^{\frac{1}{2}}$, which is also positive definite, one concludes that all the eigenvalues of ${\mathbf C}$ are positive.
\item[\rm \bf Step 3.] Applying the inequality between arithmetic and geometric means to the eigenvalues $\lambda_1, \ \ldots, \ \lambda_n$ of $\mathbf C$, one obtains:
    $$\det({\mathbf C}) = \lambda_1 \ldots \lambda_n \leq \left(\frac{\lambda_1 + \ldots + \lambda_n}{n}\right)= 1. $$
\end{enumerate}

The above reasoning is mainly based on the fact that the positivity of the spectrum of a symmetric positive definite matrix is preserved under multiplication by a positive diagonal matrix. Let us recall the following definition (see \cite{BAO}).

{\bf Definition 1.} A matrix ${\mathbf A} \in {\mathcal M}^{n \times n}$ is called {\it $D$-positive} if $\sigma({\mathbf D}{\mathbf A})$ is all positive, for every positive diagonal matrix $\mathbf D$.

The above proof of Hadamard inequality \eqref{hadeq} remains valid for the entire class of $D$-positive matrices. Besides symmetric positive definite matrices, this class includes totally positive and Kotelyansky matrices (for definitions and determinantal properties see \cite{BAO}, \cite{CARL}).

In this paper, we study the following more general class of matrices.

Consider the conic sector around the negative direction of the real axis with the apex at the origin and inner angle $2\theta$, $0 < \theta < \frac{\pi}{2}$,
\begin{equation}\label{sector}
{\mathbb C}^0_{\theta} = \{z = x+iy \in {\mathbb C}: x < 0; -x\tan\theta < y < x\tan\theta\}.
\end{equation}
For the extremal cases $\theta = 0$ and $\theta = \frac{\pi}{2}$, we shall consider the negative direction ${\mathbb R}^-$ of the real axis (excluding zero) and the open left hand side ${\mathbb C}^-:=\{z = x+iy \in {\mathbb C}: x < 0\}$ of the complex plane, respectively.

Let us recall the following definition, introduced in \cite{KUPA}.

{\bf Definition 2.} For a given value $\theta$, $0 < \theta < \frac{\pi}{2}$, we call an $n \times n$ real matrix $\mathbf A$ {\it relatively $D$-stable} if $\sigma({\mathbf D}{\mathbf A}) \subset {\mathbb C}^0_{\theta}$ for every positive diagonal matrix $\mathbf D$.

{\bf Definition 3.} For a given value $\theta$, $0 < \theta < \frac{\pi}{2}$, we call an $n \times n$ real matrix $\mathbf A$ {\it relatively additive $D$-stable} if $\sigma({\mathbf A} - {\mathbf D}) \subset {\mathbb C}^0_{\theta}$ for every nonnegative diagonal matrix $\mathbf D$.

These are generalizations of well-known definitions of $D$-stability and additive $D$-stability (see e.g. \cite{CROSS}).

{\bf Definition 2'.} An $n \times n$ real matrix $\mathbf A$ is called {\it $D$-stable} if $\sigma({\mathbf D}{\mathbf A}) \subset {\mathbb C}^-$ for every positive diagonal matrix $\mathbf D$.

{\bf Definition 3'.} An $n \times n$ real matrix $\mathbf A$ is called {\it additive $D$-stable} if $\sigma({\mathbf A} - {\mathbf D}) \subset {\mathbb C}^-$ for every nonnegative diagonal matrix $\mathbf D$.

The paper is organized as follows. In Section 2, we prove one of the main results of the paper:

\begin{theorem}\label{main} Let ${\mathbf A} = \{a_{ij}\}_{i,j = 1}^n \in {\mathcal M}^{n \times n}$. If $-{\mathbf A}$ is relatively $D$-stable for some $\theta$, $0 \leq \theta < \frac{\pi}{2}$, and, in addition, $a_{ii} > 0$ for each $i = 1, \ \ldots, \ n$, then
\begin{equation}\label{had} \det{\mathbf A} \leq \frac{1}{\cos^n \theta}\prod_{i=1}^na_{ii}. \end{equation}
\end{theorem}

The proof is based on the generalization of the inequality of arithmetic and geometric means for complex numbers.

In Section 3, we prove a weaker inequality for relatively additive $D$-stable matrices. In Section 4, we list some known classes of relatively $D$-stable matrices. Besides this, we measure the sector gap between the spectrum of a diagonally stable (see Definition 5) matrix and the imaginary axis (see Theorem \ref{dest}). In Section 5 we apply Theorem \ref{main} to obtain upper bounds for determinants of diagonally stable matrices. Section 6 provides upper bounds for determinants of diagonally dominant matrices. In Section 7, we obtain a generalization of Hadamard inequality for one more class of $D$-stable matrices.

\section{Generalized Hadamard inequality for relatively $D$-stable matrices}
First, let us consider the following generalization of the inequality of arithmetic and geometric means (see \cite{WIL}, p. 263, Theorem 1).
\begin{theorem}\label{mean} Let $z_1, \ \ldots, \ z_n \in {\mathbb C}$ and $|{\rm arg} (z_i)| \leq \theta < \frac{\pi}{2}$, $i = 1, \ \ldots, \ n$. Then
\begin{equation}\label{1} |z_1 \ldots z_n|^{\frac{1}{n}} \leq \frac{1}{n\cos \theta}|z_1 + \ldots + z_n|, \end{equation}
where the equality holds if and only if $n$ is even and $$z_1 = \ldots = z_{\frac{n}{2}} = \overline{z}_{\frac{n}{2}+1} = \ldots = \overline{z}_n = re^{i\theta}.$$
\end{theorem}
Now we can prove the following statement, which generalizes the classical Hadamard's inequality \eqref{hadeq}.

{\bf Proof of Theorem 2.} Since $a_{ii} > 0$ for every $i = 1, \ \ldots, \ n$, we can define a positive diagonal matrix ${\mathbf D}$ as follows:
$${\mathbf D} = {\rm diag}\{d_{11}, \ \ldots, \ d_{nn}\}, \qquad d_{ii} = \frac{1}{a_{ii}}, \qquad i = 1, \ \ldots, \ n. $$
Taking into account $\det{\mathbf D} = \frac{1}{a_{11} \ldots a_{nn}}$, we conclude that Inequality \eqref{had} is equivalent to the following inequality:
$$\det({\mathbf D}{\mathbf A}) \leq \frac{1}{\cos^n\theta}.$$
Denote the eigenvalues of ${\mathbf D}{\mathbf A}$ by $\widetilde{\lambda}_1, \ \ldots, \ \widetilde{\lambda}_n$. Since $-{\mathbf A}$ is relatively $D$-stable, any eigenvalue $\widetilde{\lambda}_i$ of ${\mathbf D}{\mathbf A}$ satisfies the inequality
$$|{\rm arg} (\widetilde{\lambda}_i)| \leq \theta < \frac{\pi}{2}, \qquad i = 1, \ \ldots, \ n. $$
Thus, applying Theorem \ref{mean}, we obtain the estimate:
$$ \det({\mathbf D}{\mathbf A}) = \prod_{i = 1}^n\widetilde{\lambda}_i \leq \left(\frac{1}{n \cos\theta}\left|\sum_{i = 1}^n\widetilde{\lambda}_i\right|\right)^n = \left(\frac{1}{n \cos\theta}{\rm Tr}({\mathbf D}{\mathbf A})\right)^n = \frac{1}{\cos^n\theta},$$
which implies Inequality \ref{had}.
 $\square$

The above bound for matrix determinants is applicable to any class of matrices with a sector gap between their spectra and the imaginary axis (excluding zero), which is preserved under multiplication by a diagonal matrix. Later, we shall consider some of such matrix classes.

\section{Generalized Hadamard inequality for relatively additive $D$-stable matrices}
First, let us recall the following definitions.

{\bf Definition 4.} A matrix ${\mathbf A} \in {\mathcal M}^{n \times n}$ is called a {\it $P$-matrix ($P_0$-matrix)} if all its principal minors are positive (respectively, nonnegative), i.e the inequality $A \left(\begin{array}{ccc}i_1 & \ldots & i_k \\ i_1 & \ldots & i_k \end{array}\right) > 0$ (respectively, $\geq 0$)
holds for all $(i_1, \ \ldots, \ i_k), \ 1 \leq i_1 < \ldots < i_k \leq n$, and all $k, \ 1 \leq k \leq n$.

{\bf Definition 4'.} A matrix ${\mathbf A} \in {\mathcal M}^{n \times n}$ is called a {\it $P_0^+$-matrix} if it is a $P_0$-matrix and, in addition, the sums of all principal minors of every fixed order $i$ are positive $(i = 1, \ \ldots, \ n)$.

Let us state and prove the following lemma which describes some properties of the determinant of a $P_0$-matrix. Here and later on we use the notation $[n]$ for the ordered set of indices $\{1, \ \ldots, \ n\}$.

\begin{lemma}\label{est} Let $\mathbf A$ be a $P_0$-matrix, ${\mathbf D}$ be a nonnegative diagonal matrix. Then
$$\det{\mathbf A} + \det{\mathbf D} \leq \det({\mathbf A} + {\mathbf D}).$$
\end{lemma}
{\bf Proof.} Denote the $i$th column of $\mathbf A$ by $\overline{{\mathbf a}}_i$, and the $i$th basic vector by $\overline{{\mathbf e}}_i$, $i = 1, \ldots, \ n$. Using the linearity of exterior products, we obtain the expansion
$$\det({\mathbf A} + {\mathbf D}) = (\overline{{\mathbf a}}_1 + d_{11}\overline{{\mathbf e}}_1)\wedge \ldots \wedge(\overline{{\mathbf a}}_n + d_{nn}\overline{{\mathbf e}}_n) = $$
\begin{equation}\label{!}\sum_{k=0}^n\sum_{(i_1,\ldots, i_k)\subseteq [n]}d_{i_1i_1}\ldots d_{i_ki_k}{\mathbf A}\begin{pmatrix}j_1 &  \ldots & j_{n-k} \\ j_1 &  \ldots & j_{n-k} \end{pmatrix}, \end{equation}
where the two sets $(i_1,\ldots, i_k)$, $1 \leq i_1 < \ldots < i_k \leq n$, and $(j_1,\ldots, j_{n-k})$, $1 \leq j_1 < \ldots < j_{n-k} \leq n$ are connected as follows: $$(j_1,\ldots, j_{n-k}) = [n]\setminus (i_1,\ldots, i_k).$$

Since $\mathbf A$ is a $P_0$-matrix, we have ${\mathbf A}\begin{pmatrix}j_1 &  \ldots & j_{n-k} \\ j_1 &  \ldots & j_{n-k} \end{pmatrix} \geq 0$ for all $(j_1,\ldots, j_{n-k})$, $1 \leq j_1 < \ldots < j_{n-k} \leq n$. Thus, extracting the first and last term from \eqref{!}, and taking into account $d_{ii} \geq 0$ for all $i = 1, \ \ldots, \ n$, we obtain
$$\det({\mathbf A} + {\mathbf D}) = \det{\mathbf A} + \sum_{k=1}^{n-1}\sum_{(i_1,\ldots, i_k)\subseteq [n]}d_{i_1i_1}\ldots d_{i_ki_k}{\mathbf A}\begin{pmatrix}j_1 &  \ldots & j_{n-k} \\ j_1 &  \ldots & j_{n-k} \end{pmatrix} + \det{\mathbf D}  $$ $$\geq \det{\mathbf A} + \det{\mathbf D}.$$ $\square$

Later on, we shall need the basic properties of relatively additive $D$-stable matrices, which are in fact the same as basic properties of additive $D$-stable matrices.

\begin{lemma}\label{elem}
Let ${\mathbf A} \in {\mathcal M}^{n \times n}$ be relatively additive $D$-stable for some given value $\theta$, $0 < \theta < \frac{\pi}{2}$. Then each of the following matrices is also relatively additive $D$-stable:
\begin{enumerate}
\item[\rm{(i)}] ${\mathbf A}^T$;
\item[\rm{(ii)}] ${\mathbf P}^T{\mathbf A}{\mathbf P}$, where $\mathbf P$ is a permutation matrix;
\item[\rm{(iii)}] ${\mathbf E}{\mathbf A}{\mathbf E}^{-1}$, where $\mathbf E$ is a diagonally positive matrix;
\item[\rm{(iv)}] ${\mathbf A} + {\mathbf E}$, where $\mathbf E$ is a diagonally nonnegative matrix;
\item[\rm{(v)}] $\alpha{\mathbf A}$, where $\alpha > 0$.
\end{enumerate}
\end{lemma}
{\bf Proof.} The proof copies the proof of the corresponding properties of additive $D$-stable matrices. For the proof of $(i)-(iv)$, it is enough to notice the following:
$${\mathbf A}^T - {\mathbf D} = ({\mathbf A} - {\mathbf D}^T)^T;$$
$${\mathbf P}^T{\mathbf A}{\mathbf P} - {\mathbf D} = {\mathbf P}^T({\mathbf A} - {\mathbf P}{\mathbf D}{\mathbf P}^T){\mathbf P};$$
$${\mathbf E}{\mathbf A}{\mathbf E}^{-1}- {\mathbf D} = {\mathbf E}({\mathbf A} - {\mathbf E}^{-1}{\mathbf D}{\mathbf E}){\mathbf E}^{-1}$$
$${\mathbf A} - {\mathbf E} - {\mathbf D} = {\mathbf A} - ({\mathbf E} + {\mathbf D}).$$
Taking into account that, for any nonnegative diagonal matrix $\mathbf D$, the matrices ${\mathbf D}^T$, ${\mathbf P}{\mathbf D}{\mathbf P}^T$, ${\mathbf E}^{-1}{\mathbf D}{\mathbf E}$ and ${\mathbf E} + {\mathbf D}$ are also nonnegative diagonal we complete the proof.

$(v)$. Consider an arbitrary $\lambda \in \sigma(\alpha{\mathbf A} + {\mathbf D})$. In this case, $\frac{1}{\alpha}\lambda \in \sigma({\mathbf A} + \frac{1}{\alpha}{\mathbf D})$. Since $\frac{1}{\alpha}{\mathbf D}$ is also nonnegative diagonal, we obtain $\frac{1}{\alpha}\lambda \in {\mathbb C}^0_{\theta}$ which implies $\lambda \in {\mathbb C}^0_{\theta}$.
 $\square$

For relatively additive $D$-stable matrices, let us prove the following relaxation of Inequality \eqref{had}.

\begin{theorem}\label{main2} Let ${\mathbf A} = \{a_{ij}\}_{i,j = 1}^n \in {\mathcal M}^{n \times n}$. If $-{\mathbf A}$ is relatively additive $D$-stable for some $\theta$, $0 \leq \theta < \frac{\pi}{2}$, then
\begin{equation}\label{had2} \det{\mathbf A} \leq \left(\frac{\max_i a_{ii}}{\cos \theta}\right)^n. \end{equation}
\end{theorem}
{\bf Proof.} Since $-{\mathbf A}$ is relatively additive $D$-stable for some $\theta$, $0 \leq \theta < \frac{\pi}{2}$, then $-{\mathbf A}$ is additive $D$-stable and thus $\mathbf A$ is a $P_0^+$-matrix (see \cite{CROSS}, p. 256, Corollary). So we conclude $\max_i a_{ii} > 0$. Let $\widetilde{{\mathbf A}} = \frac{1}{\max_i a_{ii}}{\mathbf A}$. By Lemma \ref{elem}, $\widetilde{{\mathbf A}}$ is also relative additive $D$-stable. Let us put
$${\mathbf D}:={\rm diag}\{1 - \widetilde{a}_{11}, \ \ldots, \  1 - \widetilde{a}_{nn}\}.$$
Since $0 \leq \frac{a_{ii}}{\max_i a_{ii}} \leq 1$, we have $0 \leq 1 - \widetilde{a}_{ii}$ for all $i = 1, \ \ldots, \ n$. Thus $\mathbf D$ is a nonnegative diagonal matrix with at least one zero entry on the principal diagonal.

Denote the eigenvalues of $\widetilde{{\mathbf A}} + {\mathbf D}$ by $\widetilde{\lambda}_1, \ \ldots, \ \widetilde{\lambda}_n$. Since $-\widetilde{{\mathbf A}}$ is relatively additive $D$-stable, any eigenvalue $\widetilde{\lambda}_i$ of $\widetilde{{\mathbf A}} + {\mathbf D}$ satisfies the inequality
$$|{\rm arg} (\widetilde{\lambda}_i)| \leq \theta < \frac{\pi}{2}, \qquad i = 1, \ \ldots, \ n. $$
Thus, applying consequently Lemma \ref{est} and Theorem \ref{mean}, we obtain the estimate:
$$\det \widetilde{{\mathbf A}}+ \det{\mathbf D} \leq \det(\widetilde{{\mathbf A}}+ {\mathbf D}) = \prod_{i = 1}^n\widetilde{\lambda}_i \leq \left(\frac{1}{n \cos\theta}\sum_{i = 1}^n\widetilde{\lambda}_i\right)^n = $$ $$ \left(\frac{1}{n \cos\theta}{\rm Tr}(\widetilde{{\mathbf A}} + {\mathbf D})\right)^n = \frac{1}{\cos^n\theta},$$
Taking into account $\det{\mathbf D} = \prod_{i = 1}^n(1 - \frac{a_{ii}}{\max_i a_{ii}})=0 $ and $\det{\widetilde{\mathbf A}} = \frac{\det{\mathbf A}}{(\max_i a_{ii})^n}$, we obtain the following inequality:
$$\frac{1}{(\max_i a_{ii})^n} \det{\mathbf A} \leq \frac{1}{\cos^n\theta},$$ which obviously implies Inequality \ref{had2}.
$\square$

\section{Classes of relatively $D$-stable and relatively additive $D$-stable matrices}
Let us consider the following four classes of matrices which will be shown to be relatively $D$-stable.
\subsection{Diagonally ${\mathbb C}^{0}_{\theta}$-stable matrices}
Let us recall the following important class of $D$-stable matrices (see \cite{CROSS}).

{\bf Definition 5.} A matrix ${\mathbf A} \in {\mathcal M}^{n \times n}$ is called {\it diagonally stable}, if there is a positive diagonal matrix ${\mathbf D}$ such that ${\mathbf W} = {\mathbf D}{\mathbf A} + {\mathbf A}^T{\mathbf D}$ is negative definite. Here, the matrix ${\mathbf D}$ is called a {\it Lyapunov scaling factor}.

By analogy with diagonally stable matrices, we introduce the following class of matrices.

{\bf Definition 6.} A matrix $\mathbf A$ is called {\it diagonally ${\mathbb C}^{0}_{\theta}$-stable}, if there exists a positive diagonal matrix $\mathbf D$ such that the matrix \begin{equation}\label{dst}{\mathbf W}({\mathbf A}, {\mathbf D}) =\begin{pmatrix} \sin(\theta){\mathbf D}{\mathbf A} & \cos(\theta){\mathbf D}{\mathbf A} \\ - \cos(\theta){\mathbf D}{\mathbf A} & \sin(\theta){\mathbf D}{\mathbf A} \\ \end{pmatrix} + \begin{pmatrix} \sin(\theta){\mathbf A}^T{\mathbf D} & -\cos(\theta){\mathbf A}^T{\mathbf D} \\  \cos(\theta){\mathbf A}^T{\mathbf D} & \sin(\theta){\mathbf A}^T{\mathbf D} \\ \end{pmatrix}\end{equation}
is negative definite.

Now we shall prove that diagonally ${\mathbb C}^{0}_{\theta}$-stable are relatively $D$-stable. For this, we need the following definition and statement.

{\bf Definition 7.} A matrix ${\mathbf A} \in {\mathcal M}^{n \times n}$ is called {\it ${\mathbb C}^{0}_{\theta}$-stable} if $\sigma({\mathbf A}) \subset {\mathbb C}^{0}_{\theta}$.

Let us recall the following particular case of the generalized Lyapunov theorem (see \cite{CGA}, p. 360, Theorem 2.2).
\begin{theorem}[Lyapunov theorem for conic regions]\label{lyap}
Given a conic region ${\mathbb C}_{\theta}^0$ with $0 < \theta < \frac{\pi}{2}$, defined by \eqref{sector}, a matrix $\mathbf A$ is ${\mathbb C}_{\theta}^0$-stable if and only if there is a symmetric positive definite matrix $\mathbf H$ such that the matrix
\begin{equation}\label{sect}{\mathbf W}({\mathbf A}, {\mathbf H}) =\begin{pmatrix} \sin(\theta){\mathbf H}{\mathbf A} & \cos(\theta){\mathbf H}{\mathbf A} \\ - \cos(\theta){\mathbf H}{\mathbf A} & \sin(\theta){\mathbf H}{\mathbf A} \\ \end{pmatrix} + \begin{pmatrix} \sin(\theta){\mathbf A}^T{\mathbf H} & -\cos(\theta){\mathbf A}^T{\mathbf H} \\  \cos(\theta){\mathbf A}^T{\mathbf H} & \sin(\theta){\mathbf A}^T{\mathbf H} \\ \end{pmatrix},\end{equation}
    is negative definite.
\end{theorem}

Now we prove the generalization of the famous result by Cross (see \cite{CROSS}, p. 254, Proposition 1), which connects the concepts of diagonal stability, multiplicative and additive $D$-stability.

\begin{theorem}\label{ldst} If $\mathbf A$ is diagonally ${\mathbb C}^{0}_{\theta}$-stable for some $\theta$, $0 < \theta < \frac{\pi}{2}$, then it is both relatively $D$-stable and relatively additive $D$-stable, for the same value of $\theta$.
\end{theorem}
{\bf Proof.} The proof copies the reasoning of \cite{CROSS} (see \cite{CROSS}, p. 254, Proposition 1). Let $\mathbf A$ be diagonally ${\mathbb C}^{0}_{\theta}$-stable, i.e. there is a positive diagonal matrix ${\mathbf D}_0$, such that the matrix ${\mathbf W}({\mathbf A}, {\mathbf D}_0)$, definined by \eqref{dst} is negative definite. First, we prove that $\mathbf A$ is relative $D$-stable, i.e. that ${\mathbf D}{\mathbf A}$ is ${\mathbb C}^{0}_{\theta}$-stable for any positive diagonal matrix $\mathbf D$. For this, let us take the matrix ${\mathbf D}_0{\mathbf D}^{-1}$, which is obviously positive diagonal, thus symmetric positive definite. Taking into account that any two diagonal matrices commute, we obtain for any positive diagonal matrix ${\mathbf D}$:
$${\mathbf W}({\mathbf D}{\mathbf A}, {\mathbf D}_0{\mathbf D}^{-1}) = \begin{pmatrix} \sin(\theta)({\mathbf D}_0{\mathbf D}^{-1})({\mathbf D}{\mathbf A}) & \cos(\theta)({\mathbf D}_0{\mathbf D}^{-1})({\mathbf D}{\mathbf A}) \\ - \cos(\theta)({\mathbf D}_0{\mathbf D}^{-1})({\mathbf D}{\mathbf A}) & \sin(\theta)({\mathbf D}_0{\mathbf D}^{-1})({\mathbf D}{\mathbf A}) \\ \end{pmatrix} + $$ $$ \begin{pmatrix} \sin(\theta)({\mathbf D}{\mathbf A})^T({\mathbf D}_0{\mathbf D}^{-1}) & -\cos(\theta)({\mathbf D}{\mathbf A})^T({\mathbf D}_0{\mathbf D}^{-1}) \\  \cos(\theta)({\mathbf D}{\mathbf A})^T({\mathbf D}_0{\mathbf D}^{-1}) & \sin(\theta)({\mathbf D}{\mathbf A})^T({\mathbf D}_0{\mathbf D}^{-1}) \\ \end{pmatrix} = $$ $$\begin{pmatrix} \sin(\theta){\mathbf D}_0{\mathbf A} & \cos(\theta){\mathbf D}_0{\mathbf A} \\ - \cos(\theta){\mathbf D}_0{\mathbf A} & \sin(\theta){\mathbf D}_0{\mathbf A} \\ \end{pmatrix} + \begin{pmatrix} \sin(\theta){\mathbf A}^T{\mathbf D}_0 & -\cos(\theta){\mathbf A}^T{\mathbf D}_0 \\  \cos(\theta){\mathbf A}^T{\mathbf D}_0 & \sin(\theta){\mathbf A}^T{\mathbf D}_0 \\ \end{pmatrix} = {\mathbf W}({\mathbf A}, {\mathbf D}_0) $$ and is negative definite. Thus by Theorem \ref{lyap}, ${\mathbf D}{\mathbf A}$ is ${\mathbb C}^{0}_{\theta}$-stable for any positive diagonal matrix $\mathbf D$.

For the proof of relative additive $D$-stability, we need to show that ${\mathbf A} - {\mathbf D}$ is relatively stable for any nonnegative matrix ${\mathbf D}$. Consider
$${\mathbf W}({\mathbf A} - {\mathbf D}, {\mathbf D}_0) = {\mathbf W}({\mathbf A}, {\mathbf D}_0) - 2\sin(\theta){\rm diag}\{{\mathbf D}{\mathbf D}_0, {\mathbf D}{\mathbf D}_0\}\prec 0,$$
since both of the matrices ${\mathbf W}({\mathbf A}, {\mathbf D}_0)$ and $- 2\sin(\theta){\rm diag}\{{\mathbf D}{\mathbf D}_0, {\mathbf D}{\mathbf D}_0\}$ are negative definite. $\square$

\subsection{Diagonally stable matrices}  Here, the {\it matrix norm} $\|{\mathbf A}\|$ of a matrix ${\mathbf A} \in {\mathcal M}^{n \times n}$ is defined as $$\|{\mathbf A}\|:=\sup_{\|{\mathbf x}\| = 1}\|{\mathbf A}{\mathbf x}\|.$$ Let us recall basic properties of the matrix norm (see, for example, \cite{BHAT2}, p. 12)
\begin{enumerate}
\item[\rm 1.] $\|{\mathbf A}^T{\mathbf A}\| = \|{\mathbf A}\|^2;$
\item[\rm 2.] $\|{\mathbf A}\| = \|{\mathbf A}^T\|;$
\item[\rm 3.] $\|{\mathbf A}{\mathbf B}\| \leq \|{\mathbf A}\|\|{\mathbf B}\|;$
\item[\rm 4.] $\|{\mathbf A}\| = s_1({\mathbf A})$, where $s_1({\mathbf A})$ is the biggest singular value of $\mathbf A$ (i.e. the eigenvalue of $({\mathbf A}^T{\mathbf A})^{\frac{1}{2}}$).
\end{enumerate}

Later on, we shall need the following result from the theory of positive definite matrices (see \cite{BHAT2}, p. 13, Proposition 1.3.2).

\begin{lemma} Let ${\mathbf A}, {\mathbf B} \in {\mathcal M}^{n \times n}$ be symmetric positive definite matrices. Then the matrix $\begin{pmatrix}{\mathbf A} & {\mathbf X} \\ {\mathbf X}^T & {\mathbf B} \end{pmatrix}$ is positive definite if and only if \begin{equation}\label{cond} {\mathbf X} = {\mathbf A}^{\frac{1}{2}}{\mathbf K}{\mathbf B}^{\frac{1}{2}} \ \mbox{for some} \ {\mathbf K}, \ \|{\mathbf K}\| < 1. \end{equation}
\end{lemma}

Basing on the above lemma, we prove the following statement.
\begin{proposition}\label{no}Let ${\mathbf A}, {\mathbf B} \in {\mathcal M}^{n \times n}$ be symmetric positive definite matrices. Then the matrix $\begin{pmatrix}{\mathbf A} & {\mathbf X} \\ {\mathbf X}^T & {\mathbf B} \end{pmatrix}$ is positive definite if the following inequality between matrix norms holds:
\begin{equation}\label{norm} \|{\mathbf X}\| < \frac{1}{\sqrt{\|{\mathbf A}^{-1}\|\|{\mathbf B}^{-1}\|}}. \end{equation}
\end{proposition}
{\bf Proof.} Extracting $\mathbf K$ from Equality \eqref{cond}, we obtain ${\mathbf K} = {\mathbf A}^{-\frac{1}{2}}{\mathbf X}{\mathbf B}^{-\frac{1}{2}}$. Thus it is enough for the proof to show that $\|{\mathbf A}^{-\frac{1}{2}}{\mathbf X}{\mathbf B}^{-\frac{1}{2}}\| < 1$. Indeed, Inequality \eqref{norm} implies
$$\sqrt{\|{\mathbf A}^{-1}\|}\|{\mathbf X}\|\sqrt{\|{\mathbf B}^{-1}\|} < 1. $$
Since both the matrices $\mathbf A$ and $\mathbf B$ are symmetric positive definite, then so are ${\mathbf A}^{-1}$ and ${\mathbf B}^{-1}$. Thus ${\mathbf A}^{-\frac{1}{2}}$ and ${\mathbf B}^{-\frac{1}{2}}$ do exist and are also symmetric positive definite. From the properties of the matix norm, we deduce that $\sqrt{\|{\mathbf A}^{-1}\|} = \sqrt{\|({\mathbf A}^{-\frac{1}{2}})^T{\mathbf A}^{-\frac{1}{2}}\|} = \|{\mathbf A}^{-\frac{1}{2}}\|$, the same reasoning shows $\sqrt{\|{\mathbf B}^{-1}\|}= \|{\mathbf B}^{-\frac{1}{2}}\|$. Thus we obtain
$$\|{\mathbf A}^{-\frac{1}{2}}\|\|{\mathbf X}\|\|{\mathbf B}^{-\frac{1}{2}}\| < 1. $$
Finally, by Property 3 of the matrix norm
$$\|{\mathbf A}^{-\frac{1}{2}}{\mathbf X}{\mathbf B}^{-\frac{1}{2}}\| \leq \|{\mathbf A}^{-\frac{1}{2}}\|\|{\mathbf X}\|\|{\mathbf B}^{-\frac{1}{2}}\|, $$
that implies $$\|{\mathbf A}^{-\frac{1}{2}}{\mathbf X}{\mathbf B}^{-\frac{1}{2}}\| < 1 .$$
 $\square$

Let a matrix ${\mathbf A} \in {\mathcal M}^{n \times n}$ be diagonally stable. Then, by Definition 5, there is a positive diagonal matrix ${\mathbf D}$ such that ${\mathbf W} = {\mathbf D}{\mathbf A} + {\mathbf A}^T{\mathbf D}$ is negative definite. The continuity reasoning shows that the Lyapunov scaling factor ${\mathbf D}$ is not unique, and we denote ${\mathcal D}({\mathbf A})$ the set of all the possible Lyapunov scaling factors of a diagonally stable matrix $\mathbf A$.

Now let us state and prove the following result connecting the concepts of diagonal stability and diagonal ${\mathbb C}^{0}_{\theta}$-stability. Here, we shall use the following notations: ${\rm Sym}({\mathbf A}) = \frac{{\mathbf A} + {\mathbf A}^T}{2}$, ${\rm Skew}({\mathbf A}) = \frac{{\mathbf A} - {\mathbf A}^T}{2}$.

\begin{theorem}\label{dest} Let a matrix ${\mathbf A} \in {\mathcal M}^{n \times n}$ be diagonally stable. Then $\mathbf A$ is diagonally ${\mathbb C}^{0}_{\theta}$-stable for any value of $\theta$ which satisfies
\begin{equation}\label{diagest}\frac{\pi}{2} > \theta > \arctan \inf_{{\mathbf D} \in {\mathcal D}({\mathbf A})}\|{\rm Skew({\mathbf D}{\mathbf A})}\|\|({\rm Sym}({\mathbf D}{\mathbf A}))^{-1}\|. \end{equation}
\end{theorem}
{\bf Proof.} Let $\mathbf A$ be diagonally stable. Then there is a set ${\mathcal D}({\mathbf A})$ of positive diagonal matrices, such that ${\mathbf W} = {\mathbf D}{\mathbf A} + {\mathbf A}^T{\mathbf D}$ is negative definite for any ${\mathbf D} \in {\mathcal D}({\mathbf A})$. For the value of $\epsilon = \tan\theta - \inf_{{\mathbf D} \in {\mathcal D}({\mathbf A})}\|{\rm Skew({\mathbf D}{\mathbf A})}\|\|({\rm Sym}({\mathbf D}{\mathbf A}))^{-1}\| > 0$, we can find a matrix ${\mathbf D}_0 \in {\mathcal D}({\mathbf A})$ such that
$$\|{\rm Skew({\mathbf D}_0{\mathbf A})}\|\|({\rm Sym}({\mathbf D}_0{\mathbf A}))^{-1}\| < $$\begin{equation}\label{!!}\inf_{{\mathbf D} \in {\mathcal D}({\mathbf A})}\|{\rm Skew({\mathbf D}{\mathbf A})}\|\|({\rm Sym}({\mathbf D}{\mathbf A}))^{-1}\| + \epsilon = \tan\theta.\end{equation}
To prove that $\mathbf A$ is diagonally ${\mathbb C}^{0}_{\theta}$-stable for any value of $\theta$ satisfying \eqref{diagest}, we consider the matrix
$${\mathbf W}({\mathbf A}, {\mathbf D}_0) =\begin{pmatrix} \sin(\theta){\mathbf D}_0{\mathbf A} & \cos(\theta){\mathbf D}_0{\mathbf A} \\ - \cos(\theta){\mathbf D}_0{\mathbf A} & \sin(\theta){\mathbf D}_0{\mathbf A} \\ \end{pmatrix} + \begin{pmatrix} \sin(\theta){\mathbf A}^T{\mathbf D}_0 & -\cos(\theta){\mathbf A}^T{\mathbf D}_0 \\  \cos(\theta){\mathbf A}^T{\mathbf D}_0 & \sin(\theta){\mathbf A}^T{\mathbf D}_0 \\ \end{pmatrix}=$$
$$\begin{pmatrix} 2\sin(\theta){\rm Sym}({\mathbf D}_0{\mathbf A}) & 2\cos(\theta){\rm Skew}({\mathbf D}_0{\mathbf A}) \\ 2\cos(\theta)({\rm Skew}({\mathbf D}_0{\mathbf A}))^T & 2\sin(\theta){\rm Sym}({\mathbf D}_0{\mathbf A}) \\ \end{pmatrix}. $$
We need to show that ${\mathbf W}({\mathbf A}, {\mathbf D}_0)$ is negative definite for any value of $\theta$ satisfying \eqref{diagest}.
Taking into account that ${\rm Sym}({\mathbf D}_0{\mathbf A})$ is negative definite for any Lyapunov scaling ${\mathbf D}_0 \in {\mathcal D}({\mathbf A})$, it will be sufficient to show that $$\cos(\theta)\|{\rm Skew}({\mathbf D}_0{\mathbf A})\| < \frac{\sin(\theta)}{\|({\rm Sym}({\mathbf D}_0{\mathbf A}))^{-1}\|}$$
and apply Proposition \ref{no}.

Indeed, since $\theta > \arctan \inf_{{\mathbf D} \in {\mathcal D}({\mathbf A})}\|{\rm Skew({\mathbf D}{\mathbf A})}\|\|({\rm Sym}({\mathbf D}{\mathbf A}))^{-1}\|$, we have $$\tan \theta > \inf_{{\mathbf D} \in {\mathcal D}({\mathbf A})}\|{\rm Skew({\mathbf D}{\mathbf A})}\|\|({\rm Sym}({\mathbf D}{\mathbf A}))^{-1}\|,$$
and by \eqref{!!},
 $$\tan\theta > \|{\rm Skew({\mathbf D}_0{\mathbf A})}\|\|({\rm Sym}({\mathbf D}_0{\mathbf A}))^{-1}\|.$$
 Then we obtain
$$\frac{\sin\theta}{\cos\theta} > \|{\rm Skew({\mathbf D}_0{\mathbf A})}\|\|({\rm Sym}({\mathbf D}_0{\mathbf A}))^{-1}\|,$$
 that implies
$$ \frac{\sin(\theta)}{\|({\rm Sym}({\mathbf D}_0{\mathbf A}))^{-1}\|} > \cos(\theta)\|{\rm Skew}({\mathbf D}_0{\mathbf A})\|.$$
Applying Proposition \ref{no}, we complete the proof.
 $\square$
 \begin{corollary}\label{corr} Let a matrix ${\mathbf A} \in {\mathcal M}^{n \times n}$ be diagonally stable. Then $\mathbf A$ is relatively $D$-stable and relatively additive $D$-stable for any value of $\theta$ which satisfies Inequality \eqref{diagest}.
 \end{corollary}
 For the proof of Corollary \ref{corr}, it is enough to apply Theorem \ref{ldst} to a diagonally ${\mathbb C}^{0}_{\theta}$-stable matrix $\mathbf A$.

 Now let us consider the following particular case of Theorem \ref{dest}.
\begin{theorem} Let a matrix ${\mathbf A} \in {\mathcal M}^{n \times n}$ be (not necessarily symmetric) positive definite (i.e. ${\rm Sym}({\mathbf A})\succ 0$).  Then $- {\mathbf A}$ is diagonally ${\mathbb C}^{0}_{\theta}$-stable for any value of $\theta$ satisfying
\begin{equation}\label{diagest1}\frac{\pi}{2} > \theta > \arctan\|{\rm Skew({\mathbf A})}\|\|({\rm Sym}({\mathbf A}))^{-1}\| = \arctan\frac{\max_{\nu \in \sigma({\rm Skew}({\mathbf A}))}|\nu|}{\min_{\mu \in \sigma({\rm Sym}({\mathbf A}))}\mu}. \end{equation}
If, in addition, $\mathbf A$ is normal, then
\begin{equation}\label{diagest2}\frac{\pi}{2} > \theta > \arctan\frac{\max_{\lambda \in \sigma({\mathbf A})}|{\rm Im}(\lambda)|}{\min_{\lambda \in \sigma({\mathbf A})}|{\rm Re}(\lambda)|}. \end{equation}
\end{theorem}
{\bf Proof.} To prove \eqref{diagest1}, first we consider a negative definite matrix as a diagonally stable one with a Lyapunov scaling factor ${\mathbf D} = {\mathbf I}$ and apply Theorem \ref{dest}. Then we recall that $({\rm Sym}({\mathbf A}))^{-1}$ and ${\rm Skew}({\mathbf A})$ are symmetric and skew--symmetric respectively, thus they are both normal. By Property (4) of the operator norm and properties of normal matrices (see e.g. \cite{HOJ}), we have:
$$\|{\rm Skew}({\mathbf A})\| = s_1({\rm Skew}({\mathbf A})) = \max_{\nu \in \sigma({\rm Skew}({\mathbf A}))}|\nu|;$$
$$\|{\rm Sym}({\mathbf A}))^{-1}\| = s_1(({\rm Sym}({\mathbf A}))^{-1}) = \max_{\mu \in \sigma(({\rm Sym}({\mathbf A}))^{-1})}|\mu|.$$
Since ${\rm Sym}({\mathbf A})$ is symmetric positive definite and so is $({\rm Sym}({\mathbf A}))^{-1}$, their spectra are positive and we have
$$\sigma(({\rm Sym}({\mathbf A}))^{-1}) = \{\frac{1}{\lambda} : \lambda \in \sigma({\rm Sym}({\mathbf A}))\}.$$
Then $$\max_{\mu \in \sigma(({\rm Sym}({\mathbf A}))^{-1})}|\mu| = \frac{1}{\min_{\mu \in \sigma({\rm Sym}({\mathbf A}))}\mu}.$$
To prove Inequality \eqref{diagest2}, it is enough to notice, that $\mathbf A$ is normal if and only if ${\rm Skew}({\mathbf A})$ and ${\rm Sym}({\mathbf A})$ do commute (see \cite{HOJ}, p. 109). Since ${\rm Sym}({\mathbf A})$ has real spectrum and the spectrum of ${\rm Skew}({\mathbf A})$ is pure imaginary, thus for each $\nu \in \sigma({\rm Skew}({\mathbf A}))$, there is $\lambda \in \sigma({\mathbf A})$ with ${\rm Im}(\lambda) = \nu$
and for each $\mu \in \sigma({\rm Sym}({\mathbf A}))$, there is $\lambda \in \sigma({\mathbf A})$ with ${\rm Re}(\lambda) = \mu$.
Thus $\max_{\nu \in \sigma({\rm Skew}({\mathbf A}))}|\nu| = \max_{\lambda \in \sigma({\mathbf A})}|{\rm Im}(\lambda)|$ and $\min_{\mu \in \sigma({\rm Sym}({\mathbf A}))}|\mu| = \min_{\lambda \in \sigma({\mathbf A})}|{\rm Re}(\lambda)|$.
 $\square$

\subsection{Matrices with $Q^2$-scalings}
As we can see from the statement of Theorem \ref{main}, we can not directly apply it to the case of $D$-stable matrices, since in this case $\theta = \frac{\pi}{2}$ and $\cos\theta = 0$. For Inequality \ref{had} to make sense, we need some gap between the spectrum of $\mathbf A$ and the imaginary axis, which will be preserved under multiplication by a positive diagonal matrix. Now we consider one more class of $D$-stable matrices which has such a gap.

First, let us recall the following result on spectra localization of $Q$-matrices (see \cite{HEK}, p.  107, Theorem 1.3, also \cite{KEL}, Corollary 1).

\begin{theorem} Let $\mathbf A$ be an $n \times n$ $Q$-matrix. Then all eigenvalues $\lambda$ of $\mathbf A$ satisfy
\begin{equation}\label{es}|\arg(\lambda)| < \pi - \frac{\pi}{n}. \end{equation}
\end{theorem}

The following theorem was proved in \cite{KU1} (see \cite{KU1}, p. 182, Theorem 1.2)

\begin{theorem}\label{dstab} Let an $n \times n$ matrix $\mathbf A$ be a $P$-matrix and $({\mathbf D}{\mathbf A})^2$ be a $Q$-matrix for every positive diagonal matrix $\mathbf D$. Then $-{\mathbf A}$ is $D$-stable.
\end{theorem}

Under the assumptions of Theorem \ref{dstab}, more precise spectra localization can be given.

\begin{theorem}\label{q}Let an $n \times n$ matrix $\mathbf A$ be a $P$-matrix and $({\mathbf D}{\mathbf A})^2$ be a $Q$-matrix for every positive diagonal matrix $\mathbf D$. Then $-{\mathbf A}$ is relatively $D$-stable with $\theta = \frac{\pi}{2} - \frac{\pi}{2n}$.
\end{theorem}
{\bf Proof.} Given an arbitrary positive diagonal matrix $\mathbf D$, let us consider $\sigma({\mathbf D}{\mathbf A})$. We need to prove that $|\arg{\lambda}| < \frac{\pi}{2} - \frac{\pi}{2n}$ for each $\lambda \in \sigma({\mathbf D}{\mathbf A})$. First, by Theorem \ref{dstab}, we obtain that $|\arg{\lambda}| < \frac{\pi}{2}$ for each $\lambda \in \sigma({\mathbf D}{\mathbf A})$. Then, applying Inequality \eqref{es} to $({\mathbf D}{\mathbf A})^2$, we obtain $|\arg(\lambda^2)| < \pi - \frac{\pi}{n}$, which, together with stability implies $|\arg{\lambda}| < \frac{\pi}{2} - \frac{\pi}{2n}$. $\square$

\subsection{Diagonally ${\mathbb C}^{0}_{\theta}$-dominant matrices}
Let us recall the following definition (see \cite{KUPA2}, Definition 6').

{\bf Definition 8.} Given a value $\theta \in (0, \frac{\pi}{2}]$, a real $n \times n$ matrix ${\mathbf A}$ is called {\it diagonally ${\mathbb C}^{0}_{\theta}$-dominant} if the following inequalities hold:
 \begin{enumerate}
\item[\rm 1.]$\sin \theta|a_{ii}| > \sum_{j\neq i}|a_{ij}| \qquad i = 1, \ \ldots, \ n.$
\item[\rm 2.] $a_{ii} < 0$, $i = 1, \ \ldots, \ n.$
\end{enumerate}

As it was proven in \cite{KUPA2} (see \cite{KUPA2}, Theorem 4), {\it a diagonally ${\mathbb C}^{0}_{\theta}$-dominant matrix is relatively $D$-stable}. To prove relative additive $D$-stability, it's enough to notice that if $\mathbf A$ is diagonally ${\mathbb C}^{0}_{\theta}$-dominant, then ${\mathbf A} - {\mathbf D}$ is also diagonally ${\mathbb C}^{0}_{\theta}$-dominant for any nonnegative diagonal matrix $\mathbf D$.

Let us also recall a well-known definition of diagonally dominant matrices.

{\bf Definition 8'.} A matrix ${\mathbf A} \in {\mathcal M}^{n \times n}({\mathbb C})$ is called {\it strictly row diagonally dominant} if the following inequalities hold:
 \begin{equation}\label{ROW}|a_{ii}| > \sum_{i\neq j}|a_{ij}| \qquad i = 1, \ \ldots, \ n.\end{equation}
 A matrix $\mathbf A$ is called {\it strictly column diagonally dominant} if ${\mathbf A}^T$ is strictly row diagonally dominant and {\it doubly diagonally dominant} if it is both row and column (strictly) diagonally dominant.
 
 It is well-known (see e.g. \cite{JOHN1}) that {\it any strictly diagonally dominant matrix with negative principal diagonal entries is $D$-stable}. Now let us prove that such a matrix is relatively $D$-stable.
 
 \begin{theorem}\label{ddom} Let ${\mathbf A} = \{a_{ij}\}_{i,j = 1}^n$ be strictly row diagonally dominant matrix with $a_{ii} < 0$, $i = 1, \ \ldots, \ n$. Denote $s_i := \frac{\sum_{i\neq j}|a_{ij}|}{|a_{ii}|}$. Then $\mathbf A$ is relatively $D$-stable for any value of $\theta$ satisfying
\begin{equation}\label{domq}
\frac{\pi}{2} > \theta > \arcsin \max_is_i.
\end{equation}
\end{theorem}
{\bf Proof}. For the proof, it is enough to show that $\mathbf A$ is diagonally ${\mathbb C}^{0}_{\theta}$-dominant for any value of $\theta$ which satisfies Inequality \eqref{domq}. Since $|a_{ii}| > \sum_{i\neq j}|a_{ij}|$ for all $i = 1, \ \ldots, \ n$, we have $$0 < \max_i s_i = \max_i \frac{\sum_{i\neq j}|a_{ij}|}{|a_{ii}|} < 1.$$
Thus the interval $(\arcsin \max_is_i , \frac{\pi}{2})$ is nonempty and any $\theta \in (\arcsin \max_is_i , \frac{\pi}{2})$ satisfies $\max_i s_i < \sin\theta < 1$. Therefore we obtain
$$\sin\theta|a_{ii}| > \max_i s_i |a_{ii}| = |a_{ii}| \max_i\frac{\sum_{i\neq j}|a_{ij}|}{|a_{ii}|} \geq |a_{ii}|\frac{\sum_{i\neq j}|a_{ij}|}{|a_{ii}|} = \sum_{i\neq j}|a_{ij}|,$$
 for all $i = 1, \ \ldots, \ n$. The above inequality shows that ${\mathbf A}$ is diagonally ${\mathbb C}^{0}_{\theta}$-dominant and, by \cite{KUPA2}, Theorem 4, relatively $D$-stable. $\square$

\section{Some bounds for the determinants of diagonally stable matrices}
Now we can state and prove the following estimate for the determinant of a diagonally stable matrix.

\begin{theorem} Given a matrix ${\mathbf A} = \{a_{ij}\}_{i,j = 1}^n$, let $-{\mathbf A}$ be diagonally stable and ${\mathcal D}({\mathbf A})$ be the set of its Lyapunov scaling factors. Then
$$\det({\mathbf A}) \leq ((\inf_{{\mathbf D} \in {\mathcal D}({\mathbf A})}\|{\rm Skew({\mathbf D}{\mathbf A})}\|\|({\rm Sym}({\mathbf D}{\mathbf A}))^{-1}\|)^2 + 1)^{\frac{n}{2}}\prod_{i=1}^na_{ii}. $$
\end{theorem}
{\bf Proof.} Applying consequently Theorems \ref{dest} and \ref{ldst}, we obtain that a diagonally stable matrix $-{\mathbf A}$ is relatively $D$-stable for any $\theta$ satisfying $$\frac{\pi}{2} > \theta > \arctan \inf_{{\mathbf D} \in {\mathcal D}({\mathbf A})}\|{\rm Skew({\mathbf D}{\mathbf A})}\|\|({\rm Sym}({\mathbf D}{\mathbf A}))^{-1}\|.$$ Turning $\theta$ to $\arctan\inf_{{\mathbf D} \in {\mathcal D}({\mathbf A})}\|{\rm Skew({\mathbf D}{\mathbf A})}\|\|({\rm Sym}({\mathbf D}{\mathbf A}))^{-1}\|$, we obtain $$\tan(\theta) \rightarrow \inf_{{\mathbf D} \in {\mathcal D}({\mathbf A})}\|{\rm Skew({\mathbf D}{\mathbf A})}\|\|({\rm Sym}({\mathbf D}{\mathbf A}))^{-1}\|.$$ Using the trigonometric identity $$\frac{1}{\cos^2(\theta)} = 1 + \tan^2(\theta),$$ we have
$$\frac{1}{\cos(\theta)} \ \mbox{turns to} \ ((\inf_{{\mathbf D} \in {\mathcal D}({\mathbf A})}\|{\rm Skew({\mathbf D}{\mathbf A})}\|\|({\rm Sym}({\mathbf D}{\mathbf A}))^{-1}\|)^2 + 1)^{\frac{1}{2}};$$
$$\frac{1}{\cos^n(\theta)} \ \mbox{turns to} \ ((\inf_{{\mathbf D} \in {\mathcal D}({\mathbf A})}\|{\rm Skew({\mathbf D}{\mathbf A})}\|\|({\rm Sym}({\mathbf D}{\mathbf A}))^{-1}\|)^2 + 1)^{\frac{n}{2}}.$$
Applying Inequality \eqref{had}, we obtain
$$\det{\mathbf A} \leq \frac{1}{\cos^n(\theta)}a_{11}\ldots a_{nn}, $$
and taking the limit in the right-hand part
$$\det{\mathbf A} \leq ((\inf_{{\mathbf D} \in {\mathcal D}({\mathbf A})}\|{\rm Skew({\mathbf D}{\mathbf A})}\|\|({\rm Sym}({\mathbf D}{\mathbf A}))^{-1}\|)^2 + 1)^{\frac{n}{2}}a_{11}\ldots a_{nn}. $$
$\square$
\begin{corollary}\label{ex} Let a matrix ${\mathbf A} \in {\mathcal M}^{n \times n}$ be (not necessarily symmetric) positive definite (i.e. ${\rm Sym}({\mathbf A}) \succ 0$).  Then
$$\det({\mathbf A}) \leq ((\|{\rm Skew({\mathbf A})}\|\|({\rm Sym}({\mathbf A}))^{-1}\|)^2 + 1)^{\frac{n}{2}}\prod_{i=1}^na_{ii}.$$
If, in addition, $\mathbf A$ is normal, then
$$\det({\mathbf A}) \leq \left(\frac{\max_{\lambda \in \sigma({\mathbf A})}{\rm Im}^2(\lambda)}{\min_{\lambda \in \sigma({\mathbf A})}{\rm Re}^2(\lambda)} + 1\right)^{\frac{n}{2}}\prod_{i=1}^na_{ii}.$$
\end{corollary}

\section{Some bounds for determinants of diagonally dominant matrices}

Here, we provide the following estimate for determinants of diagonally dominant matrices.

\begin{theorem} Let ${\mathbf A} = \{a_{ij}\}_{i,j = 1}^n$ be strictly row diagonally dominant matrix with $a_{ii} > 0$, $i = 1, \ \ldots, \ n$. Denote $s_i := \frac{\sum_{i\neq j}|a_{ij}|}{|a_{ii}|}$. Then
\begin{equation}\label{2} \det{\mathbf A} \leq \frac{1}{(1- (\max_i s_i)^2)^{\frac{n}{2}}}\prod_{i=1}^n a_{ii}. \end{equation}
\end{theorem}
{\bf Proof}. Applying Theorem \ref{ddom}, we obtain that $-{\mathbf A}$ is relatively $D$-stable for any value of $\theta$ satisfying
$$\frac{\pi}{2} > \theta > \arcsin \max_is_i.$$
 Applying Inequality \eqref{had}, we obtain
$$ \det{\mathbf A} \leq \frac{1}{(1 - \sin^2 \theta)^{\frac{n}{2}}}\prod_{i=1}^na_{ii},$$
and, by turning $\sin\theta \rightarrow \max_i s_i $, we complete the proof.
   $\square$
   
Corollary \ref{ex} implies another bound for the determinant of a doubly diagonally dominant matrix.
\begin{theorem} Let ${\mathbf A} = \{a_{ij}\}_{i,j = 1}^n$ be a doubly diagonally dominant matrix with $a_{ii} > 0$, $i = 1, \ \ldots, \ n$. Then
$$\det({\mathbf A}) \leq ((\|{\rm Skew({\mathbf A})}\|\|({\rm Sym}({\mathbf A}))^{-1}\|)^2 + 1)^{\frac{n}{2}}\prod_{i=1}^na_{ii}.$$
\end{theorem}
{\bf Proof}. For a doubly diagonally dominant matrix $\mathbf A$ with $a_{ii} > 0$, its symmetric part ${\rm Sym}({\mathbf A})$ is also doubly diagonally dominant with positive principal diagonal entries. Thus ${\rm Sym}({\mathbf A})$ is positive definite. Applying Corollary \ref{ex}, we complete the proof. $\square$
   
\section{Generalized Hadamard inequality for matrices with $Q^2$-scalings}
Now let us consider an upper bound for determinants of the following class of $D$-stable matrices.

\begin{theorem}\label{qq} Let an $n \times n$ matrix $\mathbf A$ be a $P$-matrix and $({\mathbf D}{\mathbf A})^2$ be a $Q$-matrix for every positive diagonal matrix $\mathbf D$. Then
$$\det{\mathbf A} \leq \frac{1}{\sin^n \frac{\pi}{2n}}\prod_{i=1}^na_{ii}.$$
\end{theorem}
{\bf Proof.} Applying Theorem \ref{q}, we obtain that $-{\mathbf A}$ is relatively $D$-stable with $\theta = \frac{\pi}{2} - \frac{\pi}{2n}$. Using trigonometric identities, we obtain that $\cos(\frac{\pi}{2} - \frac{\pi}{2n}) = \sin(\frac{\pi}{2n})$. Applying Inequality \eqref{had2}, we complete the proof.  $\square$

The following results for low dimensions illustrates Theorem \ref{qq}.

\begin{lemma} Let $n = 2$ and $\mathbf A$ be a $P$-matrix. Then the following conditions are equivalent.
\begin{enumerate}
\item[\rm (i)] $({\mathbf D}{\mathbf A})^2$ is a $Q$-matrix for every positive diagonal matrix $\mathbf D$;
\item[\rm (ii)] $\det ({\mathbf A}) < \frac{1}{\sin^2\frac{\pi}{4}}a_{11}a_{22} = 2a_{11}a_{22}$.
\end{enumerate}
\end{lemma}
{\bf Proof.} Consider ${\mathbf A} = \begin{pmatrix} a_{11} & a_{12} \\ a_{21} & a_{22} \end{pmatrix}$ and ${\mathbf D} = {\rm diag}\{d_{11}, d_{22}\}$. In this case,  $${\mathbf D}{\mathbf A} = \begin{pmatrix} d_{11}a_{11} & d_{11}a_{12} \\ d_{22}a_{21} & d_{22}a_{22} \end{pmatrix}, $$ $$ ({\mathbf D}{\mathbf A})^2 = \begin{pmatrix} d_{11}^2a_{11}^2 + d_{11}d_{22}a_{12}a_{21} & * \\ * & d_{22}^2a_{22}^2 + d_{11}d_{22}a_{12}a_{21} \end{pmatrix}.$$
Then $\det(({\mathbf D}{\mathbf A})^2) = d_{11}^2d_{22}^2\det({\mathbf A}) > 0$ whenever $\det({\mathbf A}) > 0$.

Thus $\rm{(i)}$ is equivalent to $${\rm Tr}({\mathbf D}{\mathbf A})^2 = d_{11}^2a_{11}^2 + d_{22}^2a_{22}^2 + 2d_{11}d_{22}a_{12}a_{21} > 0$$ for all positive $d_{11}$ and $d_{22}$.

Since $$d_{11}^2a_{11}^2 + d_{22}^2a_{22}^2 + 2d_{11}d_{22}a_{12}a_{21} + 2d_{11}d_{22}a_{11}a_{22} - 2d_{11}d_{22}a_{11}a_{22} = $$
$$(d_{11}a_{11} - d_{22}a_{22})^2 + 2 d_{11}d_{22}(a_{11}a_{22} + 2a_{12}a_{21} - a_{12}a_{21})= $$
$$(d_{11}a_{11} - d_{22}a_{22})^2 + 2 d_{11}d_{22}(\det({\mathbf A}) + 2a_{12}a_{21}) \geq 0 $$
for all $d_{11}>0, \ d_{22} > 0$ whenever $2a_{12}a_{21} >0 $. This sufficient condition is also necessary since $d_{11}a_{11} - d_{22}a_{22} = 0$ for appropriate choise of positive $d_{11}, d_{22}$. Therefore the condition $\det({\mathbf A}) + 2a_{12}a_{21} > 0$ should always hold. $\square$


\begin{thebibliography}{}
\bibitem{BAO}
Y.S. Barkovsky, T.V. Ogorodnikova, {\it On matrices with positive and simple spectra}, Izvestiya SKNC VSH Natural sciences \textbf{4} (1987), 65--70.

\bibitem{BHAT2}
R. Bhatia, {\it Positive definite matrices}, Princeton University Press, 2007.

\bibitem{CARL}
D. Carlson, {\it Weakly sign-symmetric matrices and some determinantal inequalities}, Colloquium Mathematicum, {\bf XVII} (1967), pp. 123--129.

\bibitem{CGA}
M. Chilali, P. Gahinet, P. Apkarian, {\it Robust pole placement in LMI regions},
Proceedings of the 36th Conference on Decision and Control
San Diego, USA, pp. 1291--1296 (1997).

\bibitem{CROSS}
G.W. Cross, {\it Three types of matrix stability}, Linear Algebra Appl., {\bf 20} (1978), pp. 253--263.

\bibitem{HEK}
D. Hershkowitz, N. Keller, {\it Positivity of principal minors, sign symmetry and stability}, Linear Algebra Appl. {\bf 364} (2003), pp. 105–124.

\bibitem{HOJ}
 R. Horn, C.R. Johnson, {\it Matrix analysis,} Cambridge University Press, 1990.

\bibitem{JOHN1}
C.R. Johnson, {\it Sufficient conditions for $D$-stability,} Journal of Economic Theory {\bf 9} (1974), 53-62.

\bibitem{KEL}
R.B. Kellog, {\it On complex eigenvalues of $M-$ and $P$-matrices}, Numer. Math. {\bf 19} (1972), pp. 170–175.

\bibitem{KU1}
O.Y. Kushel, {\it On a criterion of $D$-stabiity for $P$-matrices}, Special Matrices, {\bf 4} (2016), pp. 181-188.

\bibitem{KUPA}
O.Y. Kushel, R. Pavani, {\it The problem of generalized $D$-stability in unbounded LMI regions and its computational aspects,} J. Dyn. Diff. Equat., {\bf 34}, pp. 651--669 (2022).

\bibitem{KUPA2}
O.Y. Kushel, R. Pavani, {\it Generalization of the concept of diagonal dominance with applications to matrix $D$-stability,} Linear Algebra and its Applications, {\bf 630}, pp. 204-224 (2021).

\bibitem{LICHE}
W. Li, Y. Chen, {\it Some new two-sided bounds for determinants of diagonally dominant matrices,} J. Inequal Appl., {\bf 61} (2012). https://doi.org/10.1186/1029-242X-2012-61

\bibitem{MIR}
 L. Mirsky, {\it An introduction to linear algebra,} Oxford University Press, 1955.

\bibitem{QR}
J.P. Quirk, R. Ruppert, {\it Qualitative economics and the stability of equilibrium}, Rev. Econom. Studies {\bf 32}, pp. 311--326 (1965).

\bibitem{WIL}
H.S. Wilf, {\it Some applications of the inequality of arithmetic and geometric means to polynomial
equations,} Proceedings of the American Mathematical Society, {\bf 14}, pp. 263--265 (1963).
\end{thebibliography}
\end{document}